\documentclass[preprint]{elsarticle}

\usepackage{amsmath,amssymb,bm}
\numberwithin{equation}{section} 
\usepackage{amsthm} 
\newtheorem{theorem}{Theorem}

\usepackage{diagbox}
\usepackage{float} 
\usepackage{algpseudocode}
\usepackage[linesnumbered,ruled,vlined]{algorithm2e}
\usepackage[margin=0.7in]{geometry}
\usepackage{hyperref}
\hypersetup{
	colorlinks=true,
	linkcolor=green,
	filecolor=magenta,      
	urlcolor=cyan,
	pdftitle={Overleaf Example},
	pdfpagemode=FullScreen,
}
\usepackage{ulem} 

\begin{document}
	
	\begin{frontmatter} 
		
		\title{Integral Boundary Conditions in Phase Field Models}%

		\author[1]{ Xiaofeng Xu }
		
		\author[2]{Lian Zhang}
		\author[3]{Yin Shi}
		\author[3]{Long-Qing Chen}
		\author[1]{Jinchao Xu \footnote{Corresponding author: xu@math.psu.edu (Jinchao Xu)} }
		\address[1]{Department of Mathematics, Pennsylvania State University, University Park, PA, 16802, USA}
        \address[2]{Shenzhen International Center for Industrial and Applied Mathematics, Shenzhen
Research Institute of Big Data, Shenzhen, Guangdong, 518172, China}
		\address[3]{Department of Material Science and Engineering, Pennsylvania State University, University Park, PA, 16802, USA}
		
		

		\begin{abstract}
			Modeling the chemical, electric and thermal transport as well as phase transitions and the accompanying mesoscale microstructure evolution within a material in an electronic device setting involves the solution of partial differential equations often with integral boundary conditions.
			Employing the familiar Poisson equation describing the electric potential evolution in a material exhibiting insulator to metal transitions, we exploit a special property of such an integral boundary condition, and we properly formulate the variational problem and establish its well-posedness.
			We then compare our method with the commonly-used Lagrange multiplier method that can also handle such boundary conditions. 
			Numerical experiments demonstrate that our new method achieves optimal convergence rate in contrast to the conventional Lagrange multiplier method.
			Furthermore, the linear system derived from our method is symmetric positive definite, and can be efficiently solved by Conjugate Gradient method with algebraic multigrid preconditioning.

		\end{abstract}
		
	\end{frontmatter}
	
	\section{Introduction}\label{Mathematical model}

	In this paper, we introduce a type of special integral boundary condition for the Poisson equation which one may encounter in many electric circuit systems. Our focus is on the numerical implementation of such a boundary condition. The model partial differential equation (PDE) reads as 
	
		\begin{equation}
		\label{simple model}
		\left\{
		\begin{aligned}-\nabla \cdot (\sigma \nabla \Phi)&= f
			 &&\hbox{~~~~in $\Omega = (0,L)^2 $},
			\\ \Phi &=U -R\int_{\Gamma_1} \sigma \frac{\partial \Phi}{\partial n} ds &&\hbox{~~~~on $\Gamma_1 = \{ (x,y) \in \partial \Omega : x=0\}  $} ,
			\\ \sigma \frac{\partial \Phi}{\partial n}  &= 0 &&\hbox{~~~~on $\Gamma_2 = \{ (x,y) \in  \partial \Omega : y= 0 ~\text{or}~ L\}  $} ,
			\\ \Phi &= \Phi_D &&\hbox{~~~~on $\Gamma_3 = \{ (x,y) \in  \partial  \Omega : x=L\}  $} ,
		\end{aligned}\right.
	\end{equation}
	where the unknown $\Phi$ is the electric potential, $\sigma$ is the conductivity, $f$ is a spatial charge distribution, $U$ is the voltage, $R$ is the resistance, and $\Phi_D$ is the electric potential on $\Gamma_3$.

	
	\begin{figure}[!ht]
		\label{fig: simple circuit} 
		\centering
		\includegraphics[scale=0.40]{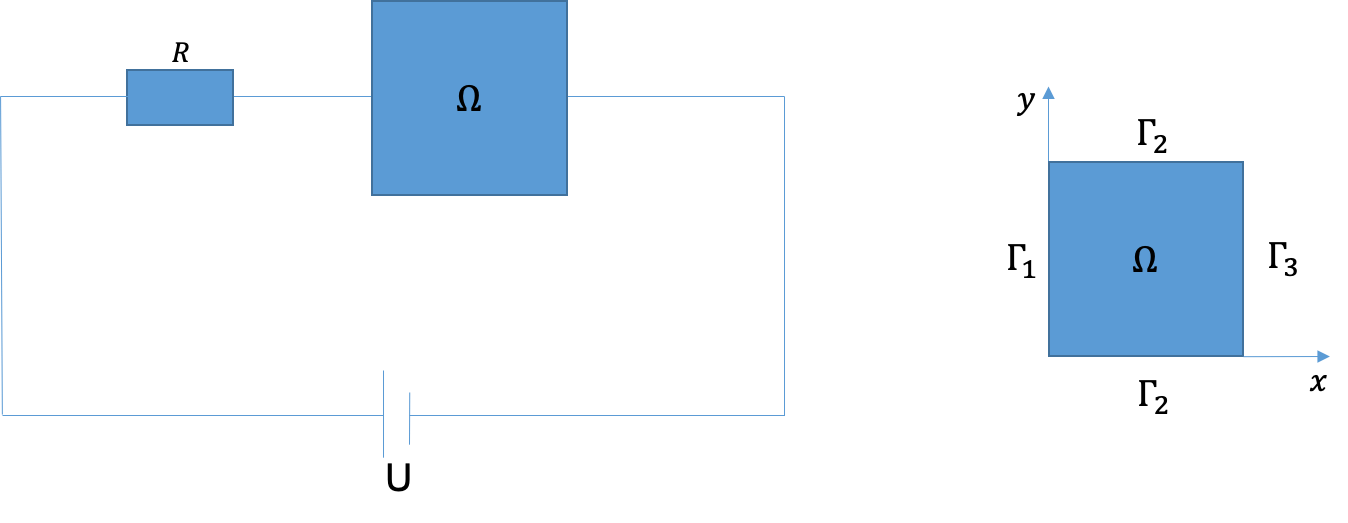}
		\caption{A simple model}
	\end{figure}

	The model PDE (\ref{simple model}) is illustrated by a typical circuit shown in Figure \ref{fig: simple circuit}, in which a material of interest is electrically excited by a direct voltage source $U$ through a series resistor $R$. 
	For simplicity, we consider a two-dimensional system with length and width both $L$ and uniform materials properties along the third dimension (out of plane). 
This PDE is frequently encountered in many simulations of mesoscale electrical systems such as rechargeable batteries~\cite{wang2020application} and resistive random access memories~\cite{zhang2020high}.
	
	We now derive the special integral boundary condition that arises on the boundary $\Gamma_1$.
	The value of $\Phi$ on $\Gamma_1$ is determined by Kirchhoff's law and Ohm's law,
	\begin{equation}
	    \Phi = U -IR,
	\end{equation}
	where $I$ is the current passing through the boundary, $R$ is the resistance of the resistor. $\Phi$ is the difference between the total voltage $U$ and the voltage drop across the resistor, $IR$.
	We denote the area and the outward unit normal vector of the boundary $\Gamma_1$ by $S=\Gamma_1\times [0,1]$ and $\bm{n}$, respectively. The current $I$ is related to the current density $\bm{J} = -\sigma \nabla\Phi$ (Ohm's law) through the following surface integral: 
	\begin{align}
		I &= -\int_S \bm{J} \cdot  \boldsymbol{n} dA \\
		& = \int_{\Gamma_1} \sigma \frac{\partial \Phi}{\partial n} ds. 
	\end{align}
	The negative sign in the first equality arises since the current flows into the material across $\Gamma_1$. Finally, we obtain the following integral boundary condition at $\Gamma_1$:   
	\begin{align} \label{integral boundary condition}
		\Phi|_{\Gamma_1} =U -R\int_{\Gamma_1} \sigma \frac{\partial \Phi}{\partial n} ds.
	\end{align}
	The potential on $\Gamma_1$ is determined by an integral of its gradient over that boundary.
	
	This type of integral boundary condition as illustrated in the problem (\ref{simple model}) is common in theoretical models of many physical processes that are controlled by external constraints. 
	For example, in the electrically induced insulator-metal transitions~\cite{Lee08Metal,Kim10Electrical,Shukla14Synchronized,Shi18Phase,Shi19Current,Kumar17Chaotic,Kumar20Third, Sood_2021}, the phase-changing material is often connected in series to a resistor. The purpose of such a resistor is to protect the material from large current damage or to control the emergence of desired phenomena such as the voltage oscillation in vanadium dioxide~\cite{Lee08Metal,Kim10Electrical,Shukla14Synchronized, PhysRevB.104.064308} and the chaotic dynamics in niobium dioxide Mott memristors~\cite{Kumar17Chaotic,Kumar20Third}.
	In simulations of these electrical devices, if not using the integral boundary condition, one must consider the system combining the material of interest and the series resistor as a whole. This results in a much more complex problem to solve, compared to the strategy to focus only on the material of interest and take account of the series resistor by the integral boundary condition. However, the integral boundary condition is a non-standard boundary condition that is not easy to deal with. Therefore, a fast and accurate method to address the integral boundary condition will greatly facilitate the simulations of various electrical devices, and is thus expected to attract broad interest of the electrical-device-modeling community.

%

To our knowledge, this type of integral boundary condition has not been analyzed theoretically or properly handled numerically in literature. 
To solve the problem by the finite element method, we need to incorporate the integral boundary condition into a weak form.
A commonly-used method is to treat the integral boundary condition as a special constraint on the boundary and enforce the constraint using the conventional Lagrange multiplier method.
Such a method often leads to a saddle point system that requires extra effort to solve with iterative solvers.
An example of treating a Dirichlet boundary condition using the Lagrange multiplier method can be found in \cite{LagmultiplierDirichlet, glowinski1995lagrange}. 
In this work, we find a new way to incorporate this special integral boundary condition directly into the weak form, leading to a symmetric positive definite linear system, which can be efficiently solved by Preconditioned Conjugate Gradient method (PCG) with an algebraic multigrid (AMG) preconditioner. We also demonstrate our proposed method is numerically superior to the existing Lagrange multiplier method in terms of accuracy.

The paper is organized as follows. In Section \ref{methods}, we formulate the variational problem using the Lagrange multiplier method, then we introduce how to incorporate the special boundary condition into the variational form naturally by using a special function space. Theoretical analysis is carried out in Section \ref{sec: analysis}, where we prove the well-posedness and the error estimates for the variational problem from our new method. 
 Numerical tests are presented in Section \ref{numericalTest}. The results show that our proposed method achieves better accuracy than the Lagrange multiplier method does, and the resulting linear system can be efficiently solved by PCG. Finally, we come to a conclusion in Section \ref{conclusion}.
	
	\section{Methods} \label{methods}
	In this section, we briefly review the conventional Lagrange multiplier method and propose our own method to solve problem (\ref{simple model}). In order to derive the variational form, we define the following function spaces:
	\begin{align}
		V_D &= \left\{ \Phi \in H^1({\Omega})~|~ \Phi|_{\Gamma_3} = \Phi_D   \right\},\\
		V &= \left\{ \Phi \in H^1({\Omega})~|~ \Phi|_{\Gamma_3}  = 0    \right\},\\
		W &= L^2(\Gamma_1),
	\end{align}
	where  $L^2(\Gamma_1) =  \left\{ \Phi : \Gamma_1 \to \mathbb{R} ~|~ \int_{\Gamma_1} \Phi^2 dx < \infty \right\}$, and $H^1({\Omega}) = \left\{ \Phi \in L^2(\Omega) ~|~ \int_{\Omega} |\nabla \Phi|^2 dx < \infty \right\}$. 
	
	For the space discretization, we use a uniform triangular mesh $\mathcal{T}^h$ of the domain $\Omega$, where $h$ denotes the mesh size of $\mathcal{T}^h$.  We let $V_{D,h}, V_{h}, W_h$ be the standard finite element spaces of $V_D, V, W$ defined as follows.
	
	\begin{align}
		V_{D,h} &= \left\{ \Phi_h \in C(\overline{\Omega}) \cap V_D ~|~ \Phi_h|_{\mathcal{K}} \in P_n(\mathcal{K}), \forall \mathcal{K} \in \mathcal{T}^h \right\},\\
		V_{h} &= \left\{ \Phi_h \in C(\overline{\Omega}) \cap V  ~|~ \Phi_h|_{\mathcal{K}} \in P_n(\mathcal{K}), \forall \mathcal{K} \in \mathcal{T}^h \right\},\\
		W_h &= \left\{ \lambda_h \in  C(\Gamma_1) ~|~ \lambda_h|_{\mathcal{K}\cap \Gamma_1} \in P_n(\mathcal{K}\cap \Gamma_1) , \forall \mathcal{K} \in \mathcal{T}^h\right\},
	\end{align}
	where $P_n(\mathcal{K})$ is the space of polynomials of degree at most $n$ on $\mathcal{K} \in \mathcal{T}^h $, with $n \geq 1$.  In our error estimates analysis and numerical tests, we use continuous piecewise linear functions for the finite element space, i.e., $n = 1$.

	\subsection{Lagrange multiplier method}
	The classical technique for \eqref{simple model} is to minimize the quadratic functional 
	\begin{equation}
	    E(\Phi) = \int_\Omega \bigg( \frac{1}{2}\sigma |\nabla \Phi|^2 - f\Phi \bigg) dx 
	\end{equation}
	over functions that satisfy the prescribed boundary conditions on $\partial \Omega = \Gamma_1\sqcup \Gamma_2 \sqcup \Gamma_3$. 
	By the Lagrange multiplier method \cite{LagmultiplierDirichlet,CalculusVariation}, the solution of equation (\ref{simple model}) becomes the stationary point of the following functional 
	\begin{equation}
		\begin{split}
			\int _\Omega \left(\frac{1}{2}\sigma |\nabla \Phi|^2 - f\Phi \right) dx - \int_{\Gamma_1}\lambda \left( U - \Phi  - R\int_{\Gamma_1} \sigma \frac{\partial \Phi}{\partial n} ds \right)ds,
		\end{split}
		\label{Lagrange multiplier functional}
	\end{equation}
	where $\Phi \in V_D$, and $\lambda \in W$ is the Lagrange multiplier.
	Based on the functional (\ref{Lagrange multiplier functional}), we introduce the following bilinear form 
	\begin{align}\label{Bilinear form 1}
		B(\Phi,\lambda;v,w) &=\int_{\Omega} \sigma \nabla \Phi \cdot \nabla v dx + \int_{\Gamma_1}\lambda v ds
		\notag \\
		& + R\int_{\Gamma_1}\lambda\bigg(\int_{\Gamma_1}\sigma \frac{\partial v}{\partial n}ds\bigg)ds+ \int_{\Gamma_1}\Phi w ds+R \int_{\Gamma_1}\bigg(\int_{\Gamma_1}\sigma \frac{\partial \Phi}{\partial n} ds\bigg)w ds ,
	\end{align}
	and a functional 
	\begin{align}
		F(v,w) =\int_{\Omega} f v dx + U\int_{\Gamma_1} w ds,
	\end{align}
	where $v \in V$, and $w \in W$.
	
	The stationary point $(\Phi $, $\lambda) \in  V_D \times W $ of (\ref{Lagrange multiplier functional}) is such that for all $v \in V$, $w \in W$,  
	\begin{align}\label{weak form Lagrange 1}
		B(\Phi,\lambda;v,w) = F(v,w).
	\end{align}
	The corresponding finite element formulation becomes the following: 
		
		Find $ \Phi_h \in V_{D,h} $, $\lambda_h \in W_h$ such that for all $ v_{h} \in V_{h}$, $w_{h} \in W_h$, 
		\begin{align}\label{discrete weak form Lagrange 1}
			B(\Phi_h,\lambda_h;v_{h},w_{h}) = F(v_{h},w_{h}).
	\end{align} 
	
    It is easy to see that the Lagrange multiplier method leads to a saddle point system.

	\subsection{New method}

		To derive the variational form for the original problem \eqref{simple model}, we make an important observation that the electric potential $\Phi$ is always a constant to be determined on the boundary $\Gamma_1$ due to the special integral boundary condition \eqref{integral boundary condition}. 
		
		Therefore, we define the following function spaces: 
		
		\begin{align}
			V^c = \left\{ \Phi \in H^1(\Omega)~|~ \Phi|_{\Gamma_1} \in P_0(\Gamma_1),  \Phi|_{\Gamma_3} = 0  \right\},
		\end{align}
		where $P_0(\Gamma_1)$ is the space of constant functions on $\Gamma_1$. 
	    This is not a standard Dirichlet boundary condition because the value of the constant on $\Gamma_1$ is not known. We use the superscript $c$ to emphasize this special property. 
		We let $V_{h}^c$ be the standard finite element subspace of $V^c$ defined as follows.
		\begin{equation}
			V_{h}^c = \{ \Phi_h \in C(\bar{\Omega}) \cap V^c ~ | ~ \Phi_h \in P_n(\mathcal{K}), \forall \mathcal{K} \in \mathcal{T}^h \}.
		\end{equation} 
		
		The basis functions for this new function space $V^c_h$ are quite easy to obtain and we discuss them in \ref{sec: appendix b: basis for new funciton space}.
		
		Without loss of generality, we may assume $\Phi_D = 0$ in the problem \eqref{simple model}. 
		Let $\Phi$ be a solution of \eqref{simple model}. For any $v \in V^c $, we have 
		\begin{equation}\label{orginal weak form 1}
			-\int_\Omega \nabla \cdot (\sigma \nabla \Phi) v dx= \int_\Omega \sigma \nabla \Phi \cdot \nabla v dx  - \int_{\Gamma_1}\sigma \frac{\partial \Phi}{\partial n} v ds =  \int_\Omega f vdx.
		\end{equation}
		Rearranging the integral boundary condition \eqref{integral boundary condition}, we have
		\begin{equation}
			\int_{\Gamma_1} \sigma \frac{\partial \Phi}{\partial n}ds = \frac{-\Phi + U}{R} = \frac{1}{|\Gamma_1|} \int_{\Gamma_1} \frac{-\Phi +U}{R} ds.
		\end{equation}
		Since the test function $v$ is a constant on $\Gamma_1$, we have 
		\begin{equation}\label{equality 1}
			\int_{\Gamma_1} \sigma \frac{\partial \Phi}{\partial n} vds = \frac{1}{|\Gamma_1|} \int_{\Gamma_1} \frac{-\Phi +U}{R}v ds.
		\end{equation}
		Combining \eqref{orginal weak form 1} with  \eqref{equality 1} gives 
		\begin{equation}
			\int_\Omega \sigma \nabla \Phi \cdot \nabla v dx  + \frac{1}{|\Gamma_1|} \int_{\Gamma_1} \frac{\Phi -U}{R} v ds = \int_\Omega f vdx.
		\end{equation}
		
		Finally, the variational problem for \eqref{simple model} becomes: find $\Phi \in V^c$ such that for all $v \in V^c$, 
		\begin{equation}\label{original weak form 2}
			a(\Phi,v) = L(v),
		\end{equation}
		where $a(\Phi,v) = \int_\Omega \sigma \nabla \Phi \cdot \nabla v dx  + \frac{1}{|\Gamma_1| R} \int_{\Gamma_1} \Phi  v ds $, and $L(v) = \int_\Omega f vdx + \frac{U}{|\Gamma_1| R}\int_{\Gamma_1} v ds$. 
		The corresponding finite element formulation becomes the following:  
		
		 Find $\Phi_h \in V_{h}^c$ such that for all $v_h \in V_{h}^c$, 
		\begin{equation}\label{discrete original weak form 2}
			a(\Phi_h,v_h) = L(v_h).
		\end{equation}
		
	\section{Theoretical analysis}\label{sec: analysis}	
		\begin{theorem}{(Well-posedness)}
			Assuming that $f \in L^2(\Omega)$, the variational problem (\ref{original weak form 2}) for the original problem \eqref{simple model} is well-posed if $\sigma$ is a positive function in $\bar{\Omega}$. 
		\end{theorem}
		
		\begin{proof}
			Note that we have $\sigma_m$, $\sigma_M > 0$ such that 
			\begin{equation}
				\sigma_m \leq \sigma(x) \leq \sigma_M ~~~~\forall x \in \bar{\Omega}.
			\end{equation}
			
			We have for $\Phi \in V^c$,
		    \begin{equation}
		        \begin{aligned}
				a(\Phi,\Phi) & = \int_\Omega \sigma \nabla |\Phi|^2 dx + \frac{1}{|\Gamma_1|R } \int_{\Gamma_1} \Phi^2 ds \\
				   & \geq \sigma_m \int_\Omega \nabla |\Phi|^2 dx \geq \alpha \| \Phi\|_{1,\Omega} ~~~~ \text{ for some $\alpha > 0$ by Poincar\'e inequality}.
		        \end{aligned}
	        \end{equation}
	        So the bilinear form $a(\cdot,\cdot)$ in \eqref{original weak form 2} is strictly positive definite. 
	        
	        Furthermore, we have for $\Phi \in V^c$,
	        \begin{equation}
		\begin{aligned}
			a(\Phi,v) & \leq \sigma_{M} \| \nabla \Phi\|_{0,\Omega} \|\nabla v \|_{0,\Omega} + \frac{1}{|\Gamma_1|R} \| \Phi\|_{0,\Gamma_1} \|v\|_{0,\Gamma_1}\\
						& \leq \gamma \| \Phi\|_{1,\Omega} \|v \|_{1,\Omega} ~~  \text{ for some $\gamma >0$ by Poincar\'e inequality and trace theorem}.
		\end{aligned}
	\end{equation} 
	So the bilinear form $a(\cdot,\cdot)$ is continuous. 
	
	Similarly, the linear functional $L(\cdot)$ in \eqref{original weak form 2} is bounded thanks to the Poincar\'e inequality and trace theorem. Finally, the well-posedness follows from the Lax-Milgram theorem \cite{bressan2013lecture}.
		\end{proof}
	
	For the error estimates, we consider the linear finite element space as an example. 
	
	
	\begin{theorem}{($H^1$ error estimate)}
	    	Let $\Phi$ and $\Phi_h$ be the solutions of the continuous equation \eqref{original weak form 2} and the discrete equation \eqref{discrete original weak form 2} respectively. If $\Phi \in H^2(\Omega)$, we have the following error estimate: 
	\begin{equation}
		| \Phi - \Phi_h |_{1,\Omega} \lesssim h \| \Phi\|_{2,\Omega}.
	\end{equation}
	Furthermore when $H^2$ regularity result holds, we have 
	\begin{equation}
	| \Phi - \Phi_h |_{1,\Omega} \lesssim h \| f\|_{0,\Omega}.
\end{equation}
	\end{theorem}
	\begin{proof}
	 			\begin{equation}
			\begin{aligned}
				| \Phi - \Phi_h|_{1,\Omega} & \leq \| \Phi - \Phi_h \|_{1,\Omega} \lesssim \|\Phi - \Phi_I \|_{1,\Omega} ~~\text{C\'ea's lemma}\\
				& \lesssim 	| \Phi - \Phi_I|_{1,\Omega} ~~~~~ \text{Poincar\'e inequality} \\
				& \lesssim h \| \Phi \|_{2,\Omega}\\
				&\lesssim h\| f\|_{0,\Omega}, ~~~ \text{$H^2$ regularity result}
			\end{aligned}
		\end{equation}
		where $\Phi_I \in V^c_{h}$ is the nodal interpolation function of $\Phi$.
	\end{proof}
	
	\begin{theorem}{($L^2$ error estimate)}
	    Let $\Phi$ and $\Phi_h$ be the solutions of the continuous equation \eqref{original weak form 2} and the discrete equation \eqref{discrete original weak form 2} respectively. 
		Suppose the $H^2$ regularity result holds, we have the following error estimate in $L^2$ norm
		\begin{equation}
			\| \Phi - \Phi_h\|_{0,\Omega} \lesssim h^2 \| \Phi\|_{2,\Omega},
		\end{equation}
	Furthermore,
	\begin{equation}
		\| \Phi - \Phi_h \|_{0,\Omega} \lesssim h^2 \| f\|_{0,\Omega},
	\end{equation}
	\end{theorem}
	\begin{proof} Considering the problem \eqref{original weak form 2} with $U = 0$ and assuming the $H^2$ regularity result holds, we can find $w \in H^2(\Omega)$ such that 
	\begin{equation}
		a(w,v) = (\Phi - \Phi_h , v)_{L^2(\Omega)}, ~~~~\text{for all $v \in V^c$ },
	\end{equation}
	and $\| w\|_{2,\Omega} \lesssim \|\Phi - \Phi_h \|_{0,\Omega}$. Letting $w_I \in V^c_{h}$ be the nodal interpolation function of $w$ and choosing $v = \Phi - \Phi_h $, we have 
	\begin{equation}
		\begin{aligned}
			\| \Phi- \Phi_h \|^2_{0,\Omega} &= a(w, \Phi -\Phi_h) \\
													& = a(w - w_I, \Phi -\Phi_h) ~~\text{orthogonality} \\
													& \lesssim \| w - w_I\|_{1,\Omega} \| \Phi -\Phi_h\|_{1,\Omega} ~~\text{continuity} \\
						&\lesssim | w -w_I|_{1,\Omega} | \Phi -\Phi_h|_{1,\Omega}\\
						& \lesssim h \| w \|_{2,\Omega} | \Phi -\Phi_h|_{1,\Omega} \\
													& \lesssim h \|\Phi- \Phi_h  \|_{0,\Omega} | \Phi -\Phi_h|_{1,\Omega}. ~~~\text{regularity}
		\end{aligned}
	\end{equation}
	It follows that 
	\begin{equation}
		\| \Phi- \Phi_h \|_{0,\Omega} \lesssim  h | \Phi -\Phi_h|_{1,\Omega} \lesssim h^2 \|\Phi \|_{2,\Omega} \lesssim h^2 \|f\|_{0,\Omega}.
	\end{equation}
	\end{proof}

	\section{Numerical tests} \label{numericalTest}
	In this section, we carry out numerical tests to compare the Lagrange multiplier method (\ref{discrete weak form Lagrange 1}) with the proposed new method (\ref{discrete original weak form 2}). 
	We implement the two methods using FEniCS \cite{FEniCSProject}. 
	Some more details of the implementation can be found in \ref{sec: appendix a}.
	
	In our numerical tests, the computational domain in both cases is $\Omega = [0,1]^2$. 
	The mesh size is denoted by $h$. For both the new method and Lagrange multiplier method, $P1$ finite elements are used.
		Let $\Phi_e$ be the exact solution and $\Phi_h$ be the numerical solution computed on $\Omega$ with mesh size $h$.
		In both tests, we look at the $L^2$ norm of the error inside the domain, $\|\Phi_e-\Phi_h\|_{0,\Omega}$, and the $H^1$ semi-norm of the error, $|\Phi_e - \Phi_h |_{1,(\Omega)}$, which is defined as $\|\nabla(\Phi_e-\Phi_h)\|_{0,\Omega}$

	We first show that in both cases, our proposed method gives better accuracy.
	We further demonstrate that our proposed method can be solved using PCG with AMG preconditioning \cite{introToAmg}.
	
	
	{\bf Test 1.} We consider a manufactured PDE in the same form as (\ref{simple model}), with $f(x,y) = -4xy + 2x$, $\sigma = 1$, ${\color{blue}U} = 1$, $R = 1$, $\Phi_D =   \frac{2}{3}y^3 - y^2 +  \frac{5}{6} $, $\Omega = [0,1]^2$. This manufactured PDE has a unique solution $\Phi_e = \frac{2}{3}xy^3 - xy^2 + \frac{5}{6}$.

	In this test, the resulting linear system is solved by the default direct solver using Sparse LU decomposition. 
	For a sequence of tests, the initial mesh size $h$ is set to $0.1$ and then refined by a factor of $2$. 
	For each mesh refinement, we report $\|\Phi_e-\Phi_h\|_{0,\Omega}$, and  $|\Phi_e - \Phi_h |_{1,\Omega}$. 
	
	The results obtained by the Lagrange multiplier method  and  our new method are shown in Tables \ref{Test1_Laarange} and \ref{Test1_direct}, respectively.
 Numerical solutions obtained by the new method achieve optimal convergence rate of $2$nd order under the $L^2$ norm inside the domain. 
	When the Lagrange multiplier method is used,  the $L^2$ error convergence rate is not optimal, although the numerical solutions still converge. 
	For $|\Phi_e - \Phi_h|_{1,\Omega}$, both methods achieve optimal convergence rate. 
	The new method is more accurate than the Lagrange multiplier method when comparing both $\|\Phi_e - \Phi_h \|_{0,\Omega}$ and $|\Phi_e - \Phi_h|_{1,\Omega}$ in this test. 
	
	%
	%
	%
	
	\begin{table}[H] 
		\centering
		\scalebox{1.0}{
			\begin{tabular}{||c|c |c|c |c ||} 
				\hline
				$h$ & $\|\Phi_e-\Phi_h\|_{0,\Omega}$ &order &$|\Phi_e-\Phi_h|_{1,\Omega}$&order \\ [0.5ex] 
				\hline\hline
				$0.1$  & $5.10\times 10^{-3}$ & * &$2.96\times 10^{-2}$&* \\
				$0.05$   & $2.47\times 10^{-3}$&1.04   &$1.49\times 10^{-2}$& 1.00	\\
				$0.025$ & $1.22\times 10^{-3}$& 1.02 &$7.43\times 10^{-3}$&	1.00\\
				$0.0125$ & $6.05\times 10^{-4}$& 1.01 &$3.72\times 10^{-3}$&1.00	\\
				[1ex] 
				\hline
			\end{tabular}
		}
		\caption{ Test 1 Error computation: Lagrange multiplier method with the default direct solver}
		\label{Test1_Laarange}
	\end{table}

	\begin{table}[H]
		\centering
		\scalebox{1.0}{
			\begin{tabular}{||c|c |c|c |c ||} 
				\hline
				$h$ & $\|\Phi_e-\Phi_h\|_{0,\Omega}$ &order &$|\Phi_e-\Phi_h|_{1,\Omega}$ & order \\ [0.5ex] 
				\hline\hline
				$0.1$  & $7.35\times 10^{-4}$ &* &$2.84\times 10^{-2} $& *\\
				$0.05$  & $1.85\times 10^{-4}$ &1.99  &$1.43\times 10^{-2}$&1.00	\\
				$0.025$ &$4.63\times 10^{-5}$&2.00	&$7.13\times 10^{-3}$&1.00	 \\
				$0.0125$&$1.16\times10^{-5}$&2.00 &$3.57\times 10^{-3}$&1.00	\\
				[1ex] 
				\hline
			\end{tabular}
		}
		\caption{Test 1 Error computation: new method with the default direct solver}
		\label{Test1_direct}
	\end{table}

	{\bf Test 2.} We consider another manufactured PDE, with $f(x,y) = (y+1)\sin(x) \cos(\pi y) + \pi \sin(x) \sin(\pi y) + \pi^2 (y+1)\sin(x)\cos(\pi y)$, $\sigma = y + 1$, ${\color{blue}U} = 1 + \frac{2}{\pi^2}$, $R = 1 $, $\Phi_D =  \sin(1) \cos(\pi y) + 1$, $\Omega = [0,1]^2$. This manufactured PDE has a unique solution $\Phi_e = \sin(x)\cos(\pi y) +1$. 
	
	The testing results for the Lagrange multiplier method and the new method are shown in Tables \ref{Test2_Lagrange} and \ref{Test2_Nitsche}, respectively. The numerical solutions obtained by the new method again achieve optimal convergence rate of $2$nd order under the $L^2$ norm inside the domain. For solutions obtained by the Lagrange multiplier method, the $L^2$ error convergence rate is still not optimal. Furthermore, the numerical results obtained by the new method is much more accurate in this test. 

	\begin{table}[H] 
		\centering
		\scalebox{1.0}{
			\begin{tabular}{||c|c |c|c |c  ||} 
				\hline
				$h$& $\|\Phi_e-\Phi_h\|_{0,\Omega}$  &order &$|\Phi_e-\Phi_h|_{1,\Omega}$ & order  \\ [0.5ex] 
				\hline\hline
				$0.1$  &  $3.53\times 10^{-2}$ & *	& $1.63\times 10^{-1}$ & *\\
				$0.05$   & $1.75\times 10^{-2}$   & 1.01	& $8.16\times 10^{-2}$ & 0.99 \\
				$0.025$ & $8.70\times 10^{-3}$   & 1.00	& $4.08\times 10^{-2}$ & 1.00\\
				$0.0125$ & $4.34\times 10^{-3}$   & 1.00	& $2.04 \times 10^{-2}$ & 1.00 \\
				[1ex] 
				\hline
			\end{tabular}
		}
		\caption{ Test 2 Error computation: Lagrange multiplier method with  the default direct solver}
		\label{Test2_Lagrange}
	\end{table}

	\begin{table}[H]
		\centering
		\scalebox{1.0}{
			\begin{tabular}{||c|c |c|c |c ||} 
				\hline
				$h$ & $\|\Phi_e-\Phi_h\|_{0,\Omega}$  &order&$|\Phi_e-\Phi_h|_{1,\Omega}$&\\ [0.5ex] 
				\hline\hline
				$0.1$  & $4.24 \times 10^{-3} $ & *	& $1.52 \times 10^{-1}$& *\\
				$0.05$  &  $1.07 \times 10^{-3}$ &  1.99&$7.60 \times 10^{-2}$ & 1.00\\
				$0.025$ &  $2.67 \times 10^{-4}$ & 2.00&$3.80 \times 10^{-2}$ & 1.00\\
				$0.0125$&$6.68 \times 10^{-5}$  & 2.00	& $1.90 \times 10^{-2}$& 1.00 \\
				[1ex] 
				\hline
			\end{tabular}
		}
		\caption{Test 2 Error computation: new method with the default direct solver}
		\label{Test2_Nitsche}
	\end{table}
	
	\textbf{Test 3.}
	We further investigate the effectiveness of PCG with AMG preconditioning for solving the manufactured PDEs using our new method in both Test 1 and Test 2. We compare CG without preconditioning to PCG, and record in Tables \ref{Test3_Iterations_pde1} and \ref{Test3_Iterations_pde2} the number of iterations needed for convergence in Tests 1 and 2, respectively.
	The stopping criterion is when the true residual norm is smaller than $10^{-7}$.

	The iterative solver CG with an AMG preconditioner can be directly applied to the linear system resulted from our new method. 
	The number of iterations is significantly reduced compared with using CG without a preconditioner. Furthermore, the number of iterations does not increase much when the mesh is refined. This indicates the effectiveness of AMG preconditioning.
	The ability of being solved by an iterative solver with AMG preconditioning is also an advantage of our new method, which can be crucially benefiting in addressing large scale systems.
	
	
	\begin{table}[H]
		\centering
		\scalebox{1.0}{
			\begin{tabular}{||c |c |c ||} 
				\hline
				$h$ &PCG & CG\\
				[0.5ex] 
				\hline\hline
				$0.1$ & 5  & 49 \\
				$0.05$&6 & 96 \\
				$0.025$&7  & 190	\\
				$0.0125$&7 &366\\
				[1ex] 
				\hline
			\end{tabular}
		}
		\caption{New method: number of iterations needed using PCG and CG in Test 1}
		\label{Test3_Iterations_pde1}
	\end{table}
	
	\begin{table}[H]
		\centering
		\scalebox{1.0}{
			\begin{tabular}{||c| c|c||} 
				\hline
				$h$ &PCG & CG\\
				[0.5ex] 
				\hline\hline
				$0.1$  &5 &53 \\
				$0.05$  &6 & 111	\\
				$0.025$ &7 & 225		\\
				$0.0125$&	7	& 456 \\
				[1ex] 
				\hline
			\end{tabular}
		}
		\caption{New method: number of iterations needed using PCG and CG in Test 2}
		\label{Test3_Iterations_pde2}
	\end{table}

	\section{Conclusion} \label{conclusion}
	We derive a special integral boundary condition for the Poisson equation for a typical electric circuit model. 
	By exploiting a special property of such an integral boundary condition and defining a new function space, we properly formulate the variational problem and establish its well-posedness. Our method results in a symmetric positive definite linear system, and we demonstrate that the linear system can be solved efficiently using PCG with an AMG preconditioner. Furthermore, numerical solutions from our proposed method are more accurate than that from the Lagrange multiplier method.
	Theoretical analysis and numerical experiments show that our proposed method achieves optimal convergence rate under both $L^2$ norm and $H^1$ semi-norm.  
	
	\section{Acknowledgement}
	This work is supported as part of the Computational Materials Sciences Program funded by the US Department of Energy, Office of Science, Basic Energy Sciences, under Award Number DE-SC0020145.
	
	
	\section*{Appendices}
	\renewcommand{\theequation}{A.\arabic{equation}}
	\appendix
	\section{Implementation of the Lagrange multiplier method in FEniCS}\label{sec: appendix a}
	The resulting variational form (\ref{weak form Lagrange 1}) can not be directly implemented using FEniCS due to the two double integrals in (\ref{Bilinear form 1}). Hence we introduce an auxiliary variable $r = \int_{\Gamma_1} \sigma \frac{\partial \Phi}{\partial n} ds $ and reformulate the problem. It then becomes seeking for the stationary point of the following functional
	\begin{align}\label{Lagrange multiplier functional fenics}
		\int _\Omega \big(\frac{1}{2}\sigma|\nabla \Phi|^2 - f\Phi \big) dx - \int_{\Gamma_1}\lambda \big( U - \Phi - R r \big)ds ,
	\end{align}
	under the constraint that
	\begin{align}\label{auxiliary Lagrange}
		r = \int_{\Gamma_1} \sigma \frac{\partial \Phi}{\partial n} ds .  
	\end{align}
	Let $V_3$ be the space of real numbers, and let $v_3 \in V_3 $ be the test function of $r$. 
	The variational form consists of two parts. From (\ref{Lagrange multiplier functional fenics}), we have 
	\begin{align} \label{weak without w}
		\int_{\Omega} (\sigma\nabla \Phi \cdot \nabla v - fv)dx + \int_{\Gamma_1} \lambda v ds - \int_{\Gamma_1}w(U-\Phi - Rr)ds  = 0
	\end{align}
	
	From (\ref{auxiliary Lagrange}), we can derive the following weak form,
	\begin{align}
		&\int_{\Gamma_1}\big(r - \int_{\Gamma_1}\sigma \frac{\partial \Phi}{\partial n}ds\big)v_3 ds = 0\\
		\iff&\int_{\Gamma_1} rv_3ds - \int_{\Gamma_1}\sigma \frac{\partial \Phi}{\partial n}ds \int_{\Gamma_1}v_3ds =0 \\ \label{auxiliary}
		\iff &\int_{\Gamma_1}\big(r - |\Gamma_1|\sigma \frac{\partial \Phi}{\partial n} \big)v_3ds = 0,
	\end{align}
	where $|\Gamma_1|$ is the length of the boundary $\Gamma_1$.
	
	Combining (\ref{weak without w}) and (\ref{auxiliary}), we obtain the variational form in the actual implementation,
	
	\begin{align}\label{weak form fenics implementation}
		\int_{\Omega} \sigma \nabla \Phi \cdot \nabla v dx + \int_{\Gamma_1} \lambda v ds + \int_{\Gamma_1}(\Phi +  Rr)wds + \int_{\Gamma_1}\big(r - |\Gamma_1|\sigma \frac{\partial \Phi}{\partial n} \big)v_3ds = \int_{\Omega} fvdx + \int_{\Gamma_1}U wds
	\end{align}
	
	\section{Basis functions for $V_{h}^c$}\label{sec: appendix b: basis for new funciton space}
	
	We consider the continuous linear finite element space as an example for illustration. 
	
	Let $I:=\{1, 2, 3,..., N\}$ be the index set for the nodes on $\mathcal{T}^h$ and $\{\varphi_1,\varphi_2, ..., \varphi_N \}$ be the corresponding basis functions of $V_{0,h}$. Let $K := \{k_1,k_2,...,k_l\} \subset I$ denote the set of indices of the nodes on the boundary $\Gamma_1$. 
	
	Note that for any $u \in V^c_{h} \subset V_{h}$, we have 
	\begin{equation}
    	\begin{aligned}
    	    	    u & = \sum_{i = 1}^{N} \alpha_i \varphi_i \\
    	    & = \alpha \sum_{k \in K}  \varphi_k + \sum_{i \in I -K} \alpha_i\varphi_i.
    	\end{aligned}
	\end{equation}
	Moreover, if $u = 0$, it is easy to see that $\alpha = 0$, $\alpha_i = 0$ for all $i \in I - K$. Hence, the basis functions for $V^c_{h}$ are given by 
	\begin{equation}
	    \{ \sum_{k \in K}\varphi_k \} \cup \{\varphi_i\}_{i \in I - K}.
	\end{equation}
	
	\section{Data Availability}
	The raw data and processed data required to reproduce these findings are available to download from the link
	\href{https://codeocean.com/capsule/7126496/tree/v1}.
	
	\bibliographystyle{model1-num-names}
	\bibliography{ref}

\end{document}